\documentstyle[11pt]{article}
\newlength{\standardunitlength}
\newtheorem{cor}{Corollary} \newtheorem{lemma}{Lemma}
\newtheorem{theorem}{Theorem} \newtheorem{prop}{Proposition}
\newenvironment{proof}{\noindent {\sc Proof:}}{$\Box$ \vspace{2 ex}}

\begin{document}

\begin{center} {\bf New Examples of Potential Theory on Bratteli
Diagrams} \end{center}

\begin{center}
By Jason Fulman
\end{center}

\begin{center}
Stanford University
\end{center}

\begin{center}
Department of Mathematics
\end{center}

\begin{center}
Building 380, MC 2125
\end{center}

\begin{center}
Stanford, CA 94305, USA
\end{center}

\begin{center}
fulman@math.stanford.edu
\end{center}

\begin{center}
http://math.stanford.edu/$\sim$ fulman
\end{center}

\begin{center}
December 17, 1999
\end{center}

\newpage

\begin{abstract} We consider potential theory on Bratteli diagrams
arising from Macdonald polynomials. The case of Hall-Littlewood
polynomials is particularly interesting; the elements of the diagram
are partitions, the branching multiplicites are integers, the
combinatorial dimensions are Green's polynomials, and the Jordan form
of a randomly chosen unipotent upper triangular matrix over a finite
field gives rise to a harmonic function. The case of Schur functions
yields natural deformations of the Young lattice and Plancharel
measure. Many harmonic functions are constructed and algorithms for
sampling from the underlying probability measures are
given. \end{abstract}

\section{Introduction}

	Potential theory on Bratteli diagrams is a beautiful subject,
with connections to probability and representation theory. The basic
set-up is as follows (for more details see Kerov's lovely article
\cite{Ke}). One starts with a Bratteli diagram; that is an oriented
graded graph $\Gamma= \cup_{n \geq 0} \Gamma_n$ such that

\begin{enumerate}
\item $\Gamma_0$ is a single vertex $\emptyset$.
\item If the starting vertex of an edge is in $\Gamma_i$, then its end vertex is in $\Gamma_{i+1}$.
\item Every vertex has at least one outgoing edge.
\item All $\Gamma_i$ are finite.
\end{enumerate}

For two vertices $\lambda, \Lambda \in \Gamma$, one writes $\lambda
\nearrow \Lambda$ if there is an edge from $\lambda$ to
$\Lambda$. Part of the underlying data is a multiplicity function
$\kappa(\lambda,\Lambda)$. Letting the weight of a path in $\Gamma$ be
the product of the multiplicities of its edges, one defines the
dimension $dim(\Lambda)$ of a vertex $\Lambda$ to be the sum of the
weights over all maximal length paths from $\emptyset$ to $\Lambda$
(this definition clearly extend to intervals). An important concept,
which we will be defined carefully in Section \ref{examples}, is the
boundary of a branching.

	Given a Bratteli diagram with a multiplicity function, one
calls a function $\phi$ {\it harmonic} if $\phi(0)=1$, $\phi(\lambda)
\geq 0$ for all $\lambda \in \Gamma$, and \[ \phi(\lambda) =
\sum_{\Lambda: \lambda \nearrow \Lambda} \kappa(\lambda,\Lambda)
\phi(\Lambda).\] An equivalent concept is that of coherent probability
distributions. Namely a set $\{M_n\}$ of probability distributions
$M_n$ on $\Gamma_n$ is called {\it coherent} if \[ M_{n-1}(\lambda) =
\sum_{\Lambda: \lambda \nearrow \Lambda} \frac{dim(\lambda)
\kappa(\lambda,\Lambda)}{dim(\Lambda)} M_{n}(\Lambda).\] The formula
allowing one to move between the definitions is $\phi(\lambda) =
\frac{M_n(\lambda)}{dim(\lambda)}$.

	One reason the set-up is interesting from the viewpoint of
probability theory is the fact that every harmonic function can be
written as a Poisson integral over the set of extreme harmonic
functions (which is often the Martin boundary). For the Pascal lattice
(vertices of $\Gamma_n$ are pairs $(k,n)$ with $k=0,1,\cdots,n$ and
$(k,n)$ is connected to $(k,n+1)$ and $(k+1,n+1)$), this fact is the
simplest instance of de Finetti's theorem. When the multiplicity
function $\kappa$ is integer valued, one can define a sequence of
algebras $A_n$ associated to the Bratteli diagram, and harmonic
functions correspond to certain characters of the inductive limit of
the algebras $A_n$.

	Several examples of the above constuctions have been examined
in detail. These include characters of the infinite symmetric group
\cite{VK22}, Kingman's branching (related to population genetics and
to a deformation of the uniform measure on the symmetric group)
\cite{Ki12}, Jack branching (which generalizes the previous two
examples and is also related to spherical functions of the infinite
hyperoctahedral Gelfand pair) \cite{KOO}, and differential posets
\cite{GK}. The paper \cite{BO} gives an update of recent developments,
and the book \cite{GDJ} contains much of interest.

	The point of this note is to provide new examples of potential
theory on Bratteli diagrams. The most interesting such example arises
from the probabilistic study of the Jordan form of a uniformly chosen
element of $T(n)$, the group of upper triangular matrices over a
finite field, with $1$'s along the main diagonal. Although
Hall-Littlewood polynomials come into play, the underlying Bratteli
diagram is different from the Hall-Littlewood branching defined in
\cite{Ke4}. In particular, the multiplicty function is integer
valued. The Bratteli diagrams examined here arise from work of Garsia
and Haiman \cite{GH} on Macdonald polynomials; to the best of our
knowledge this is the first attempt to examine them from the viewpoint
of potential theory.

	The computation of the Martin boundary for these branchings is
a hard open problem. The two main methods for computing boundaries are
the method of positive homomorphisms and the ergodic method
\cite{Ke}. It is unclear whether the Bratteli diagrams studied here
have multiplicative branching, which blocks use of the first
method. The ergodic method relies on precise estimates for ratios of
dimensions in the Bratteli diagram; at present this is blocked by the
current combinatorial intractability of Kostka-Foulkes
polynomials. Nevertheless, in one simple case, that of Schur
functions, the ergodic method does lead to a determination of the
Martin boundary.

	The sampling algorithms given here (and indeed this whole
note) were motivated by an effort of the author \cite{F2} to
understand a probabilistic growth algorithm of Borodin \cite{B} and
Kirillov \cite{Ki1} for Jordan form of uniformly chosen elements of
$T(n)$ in terms of symmetric function theory. For the reader's
benefit, we remark that the $q,t$ hook walk in \cite{GH} (defined on
Young tableaux of a given shape) is different from the probabilistic
growth algorithms given here; the paper \cite{GNW} and some of
\cite{Kerq} however are specializations of our sampling
algorithms. Section \ref{main} gives the general construction, and
Section \ref{examples} gives examples.

\section{The General Construction} \label{main}

	 To begin we introduce some notation, as on pages 2-5 of
\cite{Mac}. Let $\lambda$ be a partition of a non-negative integer $n
= \sum_i \lambda_i$ into non-negative integral parts $\lambda_1 \geq
\lambda_2 \geq \cdots \geq 0$. The notation $\lambda \vdash n$ or
$|\lambda|=n$ will mean that $\lambda$ is a partition of $n$. Let
$m_i(\lambda)$ be the number of parts of $\lambda$ of size $i$, and
let $\lambda'$ be the partition dual to $\lambda$ in the sense that
$\lambda_i' = m_i(\lambda) + m_{i+1}(\lambda) + \cdots$. Let
$n(\lambda)$ be the quantity $\sum_{i \geq 1} (i-1) \lambda_i$. It is
also useful to define the diagram associated to $\lambda$ as the set
of points $(i,j) \in Z^2$ such that $1 \leq j \leq \lambda_i$. We use
the convention that the row index $i$ increases as one goes downward
and the column index $j$ increases as one goes across. So the diagram
of the partition $(5441)$ is:

\[ \begin{array}{c c c c c}
		. & . & . & . & .  \\
		. & . & . & . &    \\
		. & . & . & . &    \\
		. & & & &  
	  \end{array} \] For $\lambda \nearrow \Lambda$, let
$R_{\Lambda / \lambda}$ (resp. $C_{\Lambda / \lambda}$) be the squares
of $\lambda$ in the same row (resp. colmun) as the square removed from
$\lambda$ to get $\Lambda$. This notation differs from that in
\cite{Mac}. Let $a_{\lambda}(s)$, $l_{\lambda}(s)$ be the number of
cells in $\lambda$ strictly to the east and south of $s$, and let
$h_{\lambda}(s)=a_{\lambda}(s)+l_{\lambda}(s)+1$. The notation $[n]$
will mean $\frac{q^n-1}{q-1}$, the $q$-analog of the number $n$. The
symbol $\psi_{\Lambda / \lambda}'$ (as on page 341 of \cite{Mac})
denotes

\[ \prod_{s \in C_{\Lambda / \lambda}}
\frac{1-q^{a_{\Lambda}(s)}t^{l_{\Lambda}(s)+1}}
{1-q^{a_{\Lambda}(s)+1} t^{l_{\Lambda}(s)}}
\frac{1-q^{a_{\lambda}(s)+1}
t^{l_{\lambda}(s)}} {1-q^{a_{\lambda}(s)}t^{l_{\lambda}(s)+1}}.\]

	For $\lambda \vdash n$, $f^{\lambda}$ will denote the
dimension of the irreducible representation of the symmetric group
$S_n$ parameterized by $\lambda$. Let $K_{\mu \Lambda}$ be the
Kostka-Foulkes polynomial, as in Section 6.8 of \cite{Mac} and let
$P_{\Lambda}(q,t)$ be Macdonald's polynomial.

{\bf Definition 1:} For $0 \leq q <1$ and $0<t<1$, the underlying Bratteli
diagram $\Gamma$ has as level $\Gamma_n$ all partitions $\lambda$ of
$n$. For $\lambda \nearrow \Lambda$, the multiplicty function is
defined as

\[ \kappa(\lambda,\Lambda) = \prod_{s \in R_{\Lambda / \lambda}}
\frac{t^{-l_{\Lambda}(s)}-q^{a_{\Lambda}(s)+1}}
{t^{-l_{\lambda}(s)}-q^{a_{\lambda}(s)+1}} \prod_{s \in C_{\Lambda /
\lambda}} \frac{q^{a_{\Lambda}(s)}-t^{-l_{\Lambda}(s)+1}}
{q^{a_{\lambda}(s)}-t^{-l_{\lambda}(s)+1}}.\] Letting $i$ be the
column number of the square removed to go from $\lambda$ to $\Lambda$,
this can be rewritten as

\[ \frac{1}{t^{\Lambda_i'-1}} \prod_{s \in R_{\Lambda / \lambda}}
\frac{1-q^{a_{\Lambda}(s)+1} t^{l_{\Lambda}(s)}}
{1-q^{a_{\lambda}(s)+1}t^{l_{\lambda}(s)}} \prod_{s \in C_{\Lambda /
\lambda}} \frac{1-q^{a_{\Lambda}(s)}t^{l_{\Lambda}(s)+1}}
{1-q^{a_{\lambda}(s)}t^{l_{\lambda}(s)+1}}.\]

	Equation I.10 of \cite{GH} proves that

\[ dim(\Lambda)=\frac{1}{t^{n(\Lambda)}} \sum_{\mu \vdash n} f^{\mu}
K_{\mu \Lambda} (q,t).\]

{\bf Definition 2:} For $0 \leq q <1, 0<t<1$ and $0 \leq
x_1,x_2,\cdots$ such that $\sum x_i=1$, define a family $\{M_n\}$ of
probability measures on partitions of size $n$ by

\begin{eqnarray*}
M_n(\Lambda) & = & \frac{(1-q)^{|\Lambda|} P_{\Lambda}(x;q,t)
\sum_{\mu \vdash n} f^{\mu} K_{\mu \Lambda}(q,t)} {\prod_{s \in
\Lambda} (1-q^{a_{\Lambda}(s)+1} t^{l_{\Lambda}(s)})}\\
& = & \frac{(1-q)^{|\Lambda|} t^{n(\Lambda)} dim(\Lambda)} {\prod_{s \in
\Lambda} (1-q^{a_{\Lambda}(s)+1} t^{l_{\Lambda}(s)})}
\end{eqnarray*}

	It will soon be verified later that the $M_n(\lambda)$ are in
fact probability measures, and also that they are coherent with
respect to the diagram of Definition 1.

	Lemma \ref{exchangeable} gives the combinatorial analog of the
probabilistic notion of exchangeability.

\begin{lemma} \label{exchangeable} Let $\gamma(0)=\emptyset \nearrow
\gamma(1) \cdots \nearrow \gamma(n) = \Lambda$ be any path in the
Bratteli diagram. Then

\[ \prod_{j=1}^n \frac{\psi'_{\gamma(j) /
\gamma(j-1)}}{\kappa(\gamma(j-1),\gamma(j))} = \frac{(1-q)^{n}
t^{n(\Lambda)}}{\prod_{s \in \Lambda} (1-q^{a_{\Lambda}(s)+1}
t^{l_{\Lambda}(s)})}.\] In paticular, the product depends on the path only
through its endpoint. \end{lemma}

\begin{proof} Suppose that $\gamma(j)$ is obtained from $\gamma(j-1)$
by adding to column $i$. Writing everything out, one sees that

\begin{eqnarray*}
\frac{\psi'_{\gamma(j) / \gamma(j-1)}}
{\kappa(\gamma(j-1),\gamma(j))} & = & t^{\gamma(j)_i'-1} \prod_{s \in
C_{\Lambda / \lambda} \cup R_{\Lambda / \lambda}}
\frac{1-q^{a_{\lambda}(s)+1} t^{l_{\lambda}(s)}}{1-q^{a_{\Lambda}(s)+1}
t^{l_{\Lambda}(s)}}\\
& = & t^{\gamma(j)_i'-1} (1-q) \frac{\prod_{s \in \lambda} 1-q^{a_{\lambda}(s)+1} t^{l_{\lambda}(s)}} {\prod_{s \in \Lambda} 1-q^{a_{\Lambda}(s)+1} t^{l_{\Lambda}(s)}}.
\end{eqnarray*}
Using the fact that $n(\Lambda)=\sum_{j=1}^n (\gamma(j)_i'-1)$
and multiplying terms, the result follows. \end{proof}

	Theorem \ref{coherent} proves that the family $\{M_n\}$
satisfies the coherence equation. It will then be seen that the
$\{M_n\}$ are indeed probability measures.
	
\begin{theorem} \label{coherent} For any $0 \leq q<1$ and $0<t<1$ and
$0 \leq x_1,x_2,\cdots$ satisfying $\sum x_i=1$, the set $\{M_n\}$ satisfy
the equation \[ M_{n-1}(\lambda) = \sum_{\Lambda: \lambda \nearrow
\Lambda} \frac{dim(\lambda)}{dim(\Lambda)} \kappa(\lambda,\Lambda)
M_{n}(\Lambda).\] \end{theorem}

\begin{proof} By Lemma \ref{exchangeable} and the definition of
$M_n(\lambda)$, any path from $\emptyset$ to $\Lambda$ yields the
equality

\[ M_n(\Lambda) = P_{\Lambda} dim(\Lambda) \prod_{j=1}^{n}
\left(\frac{\psi'_{\gamma(j) /
\gamma(j-1)}}{\kappa(\gamma(j-1),\gamma(j)} \right).\] In the
following equations, paths from $\emptyset$ to some $\Lambda$ such
that $\lambda \nearrow \Lambda$ are chosen so as to first go to
$\lambda$ (in a way independent of $\Lambda$) and then go to
$\Lambda$. Consequently,

\begin{eqnarray*}
& & \sum_{\Lambda: \lambda \nearrow \Lambda}
\frac{dim(\lambda)}{dim(\Lambda)} \kappa(\lambda,\Lambda) M_{n}(\Lambda)\\
& = & \sum_{\Lambda: \lambda \nearrow \Lambda} dim(\lambda)
\kappa(\lambda,\Lambda) P_{\Lambda} \prod_{j=1}^{n-1} \left(\frac{\psi'_{\gamma(j) /
\gamma(j-1)}}{\kappa(\gamma(j-1),\gamma(j)} \right) \frac{\psi'_{\Lambda / \lambda}}{\kappa(\lambda,\Lambda)}\\
& = & dim(\lambda) P_{\lambda} \prod_{j=1}^{n-1}
\left(\frac{\psi'_{\gamma(j) /
\gamma(j-1)}}{\kappa(\gamma(j-1),\gamma(j))} \right)
\sum_{\Lambda: \lambda \nearrow \Lambda}
\frac{P_{\Lambda}\psi'_{\Lambda / \lambda}}{P_{\lambda}
}\\
& = & M_{n-1}(\lambda) \sum_{\Lambda: \lambda \nearrow \Lambda}
\frac{P_{\Lambda} \psi'_{\Lambda / \lambda}}{P_{\lambda}}\\
& = & M_{n-1}(\lambda)
\end{eqnarray*}
Since $\sum x_i=1$, the final equality is simply equation 6.24
on page 340 of \cite{Mac} with $r=1$ (a Pieri rule).
\end{proof}

	Corollary \ref{measure} shows that the $\{M_n\}$ are indeed
probability measures.

\begin{cor} \label{measure} The $\{M_n\}$ of Definition 2 are
probability measures. \end{cor}

\begin{proof} The second expression for $\kappa(\lambda,\Lambda)$
implies that $dim(\Lambda) \geq 0$ for all $\Lambda$. The fact that
$P_{\Lambda} \geq 0$ follows from the hypotheses on $q,t$ and the
$x$'s, together with the skew-expansion rule (equation 7.9' on page
345 of \cite{Mac}) for Macdonald polynomials. Thus $M_n(\lambda) \geq
0$ for all $\lambda$. From the definition of $M_1$ it is a probability measure. For
larger $M_n$ this follows from induction and the equation

\begin{eqnarray*}
1 & = & \sum_{\lambda \vdash n-1} M_n(\lambda)\\
& = & \sum_{\lambda \vdash n-1} \sum_{\Lambda: \lambda \nearrow \Lambda} \frac{dim(\lambda)}{dim(\Lambda)} \kappa(\lambda,\Lambda) M_n(\Lambda)\\
& = & \sum_{\Lambda \vdash n} \sum_{\lambda: \lambda \nearrow \Lambda}  \frac{dim(\lambda)}{dim(\Lambda)} \kappa(\lambda,\Lambda) M_n(\Lambda)\\
& = & \sum_{\Lambda \vdash n} M_n(\Lambda).
\end{eqnarray*}
\end{proof}

	As a consequence of the fact that the $M_n$ are coherent, we
obtain for free a method of sampling from them. This principle is
implicit in the literature (e.g. page 144 of \cite{Ke}), but there is
a surprising simplication which occurs in our examples.

\begin{prop} \label{sample} Starting from $\emptyset$, at each stage
move to a larger partition according to the rule that the chance of
going from $\lambda$ to $\Lambda$ is $\frac{P_{\Lambda} \psi'_{\Lambda
/ \lambda}}{P_{\lambda}}$. Then after $n$ steps the probability of
being at the partition $\Lambda$ is $M_n{\Lambda}$. \end{prop}

\begin{proof} In general the transition probabilities from $\lambda$
to $\Lambda$ to sample from a coherent family $\{M_n\}$ is
$\frac{\kappa(\lambda,\Lambda) M_n(\Lambda) dim(\lambda)}
{M_{n-1}(\lambda) dim(\Lambda)}$. These sum to $1$ by the definition
of coherence, and sample from $M_n$ because

\begin{eqnarray*}
& & \sum_{\gamma: \gamma(0)=\emptyset \nearrow \cdots \nearrow
\gamma(n)=\Lambda} \prod_{j=1}^n \frac{\kappa(\gamma(j-1),\gamma(j)) M_j(\gamma(j))
dim(\gamma(j-1))} {M_{j-1}(\gamma(j-1)) dim(\gamma(j))}\\
& = & M_n(\Lambda) \sum_{\gamma: \gamma(0)=\emptyset \nearrow \cdots \nearrow
\gamma(n)=\Lambda} \frac{\prod_{j=1}^n \kappa(\gamma(j-1),\gamma(j))}{dim(\Lambda)}\\
& = & M_n(\Lambda).
\end{eqnarray*} This principle together with the formula for $M_n(\Lambda)$
inside the proof of Theorem \ref{coherent}, imply the proposition. \end{proof}

	As will be seen in Section \ref{examples}, in special cases
the algorithm of Proposition \ref{sample} yields known
results. Curiously, the transition probabilities of Proposition
\ref{sample} are exactly those on page 585 of \cite{F1}, if one
conditions on each coin coming up heads once. The motivating example
there was the probabilistic study of the $z-1$ part of the Jordan form
of a random element of $GL(n,q)$.

\section{Examples} \label{examples}

	This section gives some examples of the constructions in the
previous section. Before doing so, we define the Martin boundary
$\Delta$ and Poisson kernel $\Phi:\Gamma \times \Delta \mapsto R$ of a
branching as in \cite{Ke}, which the reader should consult for a
fuller treatment. One requires that $\Delta$ is a compact topological
space and that there is a map $i:\Gamma \mapsto \Delta$ such that

\begin{enumerate}
\item For every $\omega \in \Delta$ the function
$\phi_{\Lambda}(\omega)=\Phi(\Lambda,\omega)$ is harmonic with respect
to the branching.

\item The functions $\Phi_{\Lambda}(\omega)$ are continuous and span a
dense linear subspace in the space of continuous functions on
$\Delta$.

\item For every $\omega \in \Delta$, the measures $i(dim(\Lambda)
\Phi(\lambda;\omega))$ converge weakly as $n \rightarrow \infty$ to
the point mass $\delta_{\omega}$ at $\omega$.
\end{enumerate}

	In the case of the Young lattice, the boundary $\Delta$ is the
space of pairs $(\alpha;\beta)$ such that $\alpha_1 \geq \alpha_2 \geq
\cdots \geq 0$, $\beta_1 \geq \beta_2 \geq \cdots \geq 0$ and
$\sum_{i=1}^{\infty} \alpha_i + \sum_{i=1}^{\infty} \beta_i \leq
1$. The map $i$ send a partition $\Lambda$ to
$(\frac{f_1}{n},\frac{f_2}{n}, \cdots;\frac{g_1}{n},
\frac{g_2}{n},\cdots)$ where $f_i=\Lambda_i-i+\frac{1}{2}$ and
$g_i=\Lambda_i'-i+\frac{1}{2}$. The Poisson kernel
$\Phi(\Lambda;\alpha,\beta)$ is $s_{\Lambda}(\alpha;\beta;\gamma)$
where the $s_{\Lambda}(\alpha;\beta;\gamma)$ are the extended Schur
functions defined for instance on page 147 of \cite{Ke}.

\begin{enumerate}

\item {\bf Upper triangular matrices}

	Suppose that $q=0$ and $t=\frac{1}{q}$, where this second $q$
is the size of a finite field. Further, set $x_i=\frac{1}
{q^{i-1}}-\frac{1}{q^i}$. Several simplifications take place. First,
the multiplicities have a simple description; letting $i$ be the
column to which one adds in order to go from $\lambda$ to $\Lambda$,
it follows that $\kappa(\lambda,\Lambda)=q^{\lambda_i'}+
q^{\lambda_i'-1}+ \cdots+ q^{\lambda_{i+1}'}$. This is always
integral. Second, $dim(\Lambda)$ reduces to a Green's polynomial
$Q^{\Lambda}(q)=Q^{\Lambda}_{(1^n)}(q)$ as in Section 3.7 of
\cite{Mac}. These polynomials are important in the representation
theory of the finite general linear groups.

	The third and fourth simplifications are significant enough to
be stated as propositions.

\begin{prop} \label{simplify3} With the above specializations,
$M_n(\lambda)$ is the probability that a uniformly chosen element of
$T_n(\lambda)$ has Jordan form of shape $\lambda$. \end{prop}

\begin{proof} This follows by comparison with the formula in Theorem 1
of \cite{F2}. \end{proof}

\begin{cor} (\cite{B},\cite{Ki1}) The Jordan form of a uniformly
chosen element of $T_n(\lambda)$ can be sampled from by stopping the
following procedure after $n$ steps:

	Starting with the empty partition, at each step transition
from a partition $\lambda$ to a partition $\Lambda$ by adding a dot to
column $i$ chosen according to the rules

\begin{itemize} \item $ i=1$ with probability
$\frac{1}{q^{\lambda_1'}}$ \item $i=j>1$ with probability
$\frac{1}{q^{\lambda_j'}}-\frac{1}{q^{\lambda_{j-1}'}}$ 
\end{itemize}
\end{cor}

\begin{proof} This follows easily from the following five ingredients:
Proposition \ref{sample}, Proposition \ref{simplify3}, homogeneity of
$P_{\Lambda}$ (which implies that $P_{\Lambda}(x)= (1-\frac{1}{q})
P_{\Lambda} (1,\frac{1}{q}, \frac{1}{q^2},\cdots;0,\frac{1}{q})$),
Macdonald's principal specialization formula (page 337 of \cite{Mac}),
and a piece of paper. \end{proof}

	Note that Borodin \cite{B}, has shown that the asymptotic
Jordan form of a random element of $T(n)$ has the following shape: the
longest block has size $(1-\frac{1}{p})n$, the second block has size
$(\frac{1}{p}-\frac{1}{p^2})n$, etc. This suggests to us that the
harmonic function $\frac{M_n(\Lambda)}{dim(\Lambda)}$ is extremal. Is
it extremal for other $x_i$ such that $\sum x_i=1$? This brings us to
the following

{\bf Problem:} Find the Martin boundary of the Bratteli diagram in
this example.

\item {\bf Schur functions}

	A second example of interest occurs when $q=t<1$. Letting $i$
be the column to which one adds in order to go from $\lambda$ to
$\Lambda$, it is not hard to rewrite $\kappa(\lambda,\Lambda)$ as

\[ \frac{1}{q^{\lambda_i'}} \frac{\prod_{s \in \lambda}
[h_{\Lambda}(s)]}{\prod_{s \in \lambda} [h_{\lambda}(s)]}.\] One
checks (using the fact that $f^{\Lambda}$ is the number of paths in
the Young lattice from $\emptyset$ to $\Lambda$ and that the product
of multiplicities is path independent) that the dimension also has a
nice simplification, namely $dim(\Lambda) = \frac{f^{\Lambda} \prod_{s
\in \Lambda} h[s]}{q^{n(\Lambda)}}$.

	The measure $M_n(\Lambda)$ reduces to $s_{\Lambda}
f^{\Lambda}$, where $s_{\Lambda}$ is a Schur function. Setting
$x_1=\cdots=x_n=\frac{1}{n}$ and letting $n \rightarrow \infty$, one
obtains Plancharel measure, which is important in representation
theory and random matrix theory. A method for sampling from it was
found in \cite{GNW} (see \cite{Kerq} for extensions).

	Letting $x_1=\cdots=x_n$ satisfy $\sum x_i=1$ (all other
$x_j=0$) gives a natural deformation of Plancharel measure, studied
for instance by \cite{TW}. Stanley \cite{S} shows that this measure on
partitions also arises by applying the RSK algorithm to a random
permutation distributed after a biased riffle shuffle. Since the
quantities in Proposition \ref{sample} have simple expressions under
this specialization (e.g. page 45 of \cite{Mac}), the sampling
algorithm is useful.

\begin{theorem} The Martin boundary in these examples is the same as
for the Young lattice. \end{theorem}

\begin{proof} We use the ergodic method (Section 8 of \cite{Ke}). The
map $i$ is the same as for the Young lattice and the Poisson kernel is
defined as $\frac{q^{n(\Lambda)} s_{\lambda}(\alpha,\beta;\gamma)}
{\prod_{s \in \Lambda} [h_{\Lambda}(s)]}$. To see that the first
condition of a boundary is met, recall that the function
$\phi_t(\Lambda)=lim_{n \rightarrow \infty}
\frac{dim(\Lambda,\nu_n)}{dim(\emptyset,\nu_n)}$ is harmonic if it
exists (here $\nu_n$ is a sequence of vertices of a path with each
$\nu_n \in \Gamma_n$). In fact it is true that

\[ \frac{dim(\Lambda,\nu)}{dim(\emptyset,\nu)} =
\frac{q^{n(\Lambda)}}{\prod_{s \in \Lambda} [h_{\Lambda}]}
\frac{dim^*(\nu-\Lambda)}{dim^*(\nu)},\] where $dim^*$ denotes
dimension in the Young lattice and $dim^*(\nu-\Lambda)$ is the number
of paths in the Young lattice from $\Lambda$ to $\nu$. This follows
from the observation that the product of the mulitplicities $\kappa$
along a path in the Bratteli diagram of this example depends only on
the endpoints of the path. The fact that the second condition of a
boundary is met follows from a generating function argument showing
that the $\Phi_{\Lambda}(\omega)$ separate points of the boundary and
the Stone-Wierstrass theorem. The third condition of a boundary
amounts to exactly the same condition as for the Young lattice, and
thus holds. \end{proof}

\item {\bf Jack symmetric functions}

	A third example of interest occurs by setting $q=\frac{1}
{t^{\theta}}$, and taking the limit as $t \rightarrow 1$. The
multiplicity function $\kappa(\lambda,\Lambda)$ takes the form

\[ \prod_{s \in R_{\Lambda / \lambda}} \frac{a_{\Lambda}(s)+1+\theta
l_{\Lambda}(s)} {a_{\lambda}(s)+1+\theta l_{\lambda}(s)} \prod_{s \in
C_{\Lambda / \lambda}} \frac{a_{\Lambda}(s)+\theta (l_{\Lambda}(s)+1)}
{a_{\lambda}(s)+\theta (l_{\lambda}(s)+1)}.\] We do not know of a
simple general expression for $dim(\Lambda)$.

	Setting $\theta=0$ (Kingman branching) leads to trouble with
our formulation as it amount to setting $q=1$ before taking the limit
$t \rightarrow 1$. Setting $\theta=1$ (Schur functions) has already
been considered. The case of zonal functions
(i.e. $\theta=\frac{1}{2}$) merits further investigation.

\end{enumerate}

\section{Acknowledgments} This research was supported by an NSF
Postdoctoral Fellowship.

\end{document}